\documentstyle{article}

\newtheorem{thm}{Theorem}[section]

\newcommand{\be}{\begin{equation}}
\newcommand{\ee}{\end{equation}}
\newcommand{\ben}{\begin{enumerate}}
\newcommand{\een}{\end{enumerate}}
\newcommand{\pa}{{\partial}}

\newcommand{\e}{{\epsilon}}

\title{\Huge Geometric Meanings of Curvatures in Finsler Geometry\footnote{This article is  for the 20th Winter School on Geometry and Physics at Srni in Czech Republic. It was written 
 during my visit at Institute of Mathematics and Informatics (IMI) at  University of Debrecen in Hungary. The author would like to thank Dr. S. B\'{a}cs\'{o} and Dr. L. Kozma for their great  help and hospitality. }}
\author{Zhongmin Shen}
\date{\small Dept. of Math., IUPUI, 402 N.Blackford Street, Indianapolis, IN 46202-3216, USA}

\begin{document}
\maketitle
\section{Introduction}

In Finsler geometry, we use calculus to study the geometry of regular inner metric spaces.
In this note  I will briefly discuss various curvatures and their geometric meanings
 from the metric geometry  point of view, without going into  
 the forest of tensors.

 A metric $d$ on a topological space $M$ is a function
on $ M \times M$ with the following properties
\ben
\item[(D1)] $d(p, q) \geq 0$ and equality holds only when $p=q$;
\item[(D2)] $d(p, q) \leq d(p, r) + d(r, q)$.
\een
For a Lipschitz curve $c:[a, b]\to (M, d)$, define the {\it dilation}
of $c$ at $t\in [a, b]$ by

\[ {\rm dil}_t (c) :=\limsup_{\e \to 0^+}
\sup_{-\e + t < t_1 < t_2 < t+\e} { d(c(t_1), c(t_2))\over t_2 -t_1 }.\]
We obtain  a length structure on $M$ defined by
$$ \ell_d(c) := \int_a^b {\rm dil}_t (c) dt.$$
$d$ is said to be {\it inner} if 
$$ d(p, q) = \inf_c \ell_d(c),$$
where the infimum is taken over all Lipschitz curves $c$ from $p$ to $q$. 
Traditionally, we impose the following reversibility condition on $d$
\ben
\item[(D3)] $d(p, q)=d(q, p)$.
\een
But this reversibility condition is so restrictive   that
it eliminates lots of interesting metric structures, such as the Funk metric below.

 Let $\Omega$ be a strongly convex bounded domain in ${\rm R}^n$. For 
$p, q\in \Omega$, let $\ell_{pq}$ denote the ray issuing from $p$ to $q$
passing through $q$. Define
\be
 d_f(p, q) := \ln { |z-p|\over |z - q| },\label{Funk1}
\ee
where $z\in \partial \Omega$ is the intersection point of $\ell_{pq}$
with $\partial \Omega$. 
Then $d_f$ is an inner metric on $\Omega$, which is called the
{\it Funk metric} \cite{Funk}. The Funk metric is not reversible. Set 
\be
 d_h(p, q) := {1\over 2} \Big ( d_f(p, q) + d_f(q, p) \Big ) , \ \ \ \ \ \
p, q\in \Omega . \label{Hilbert1}
\ee
We obtain a reversible inner metric which is called the {\it Hilbert metric}. 
There are many other interesting inner metrics which are not 
Riemannian. 

\bigskip

An inner metric $d$ on a manifold $M$ is said to be {\it regular}
if there is a nonnegative function $F$ on $ TM$ such that

\ben 
\item[(F0)] for any $C^1$ curve $c: [a, b]\to M$,
$
 {\rm dil}_t (c) = F( \dot{c}(t)), \ \ \  a \leq t \leq b$;
\item[(F1)] $F$ is $C^{\infty}$ on $TM\setminus \{0\}$;
\item[(F2)] For each $x\in M$,
$F_x:= F|_{T_xM}$ is a Minkowski functional on $T_xM$, i.e.,
\ben
\item[(F2a)] $F_x(\lambda y) =\lambda F_x(y)$, $\ \forall\lambda >0$, $y\in T_xM$;
\item[(F2b)] for each $y\in T_xM\setminus \{0\}$, the induced symmetric
bilinear  form  $g_y$ on $T_xM$ is an inner product,
where
\be
 g_y (u, v) := {1\over 2} {\pa^2 \over \pa s \pa t} \Big [ F^2 (y+su+tv) \Big ]|_{s=t=0},
\ \ \ \ \ \ u,v \in T_xM.\label{gy}
\ee
\een

\een

 A {\it Finsler metric} on a manifold $M$ is a nonnegative function
$F$ on $TM$ which 
satisfies (F1) and (F2).

The Funk metric $d_f$ in (\ref{Funk1})  is regular
and the induced Finsler metric $F_f$  is determined by
\be
 x+ {y \over F_f(y)} \in \pa \Omega, \ \ \ \ \ \ y\in T_x\Omega. \label{Funk2}
\ee
The Hilbert metric $d_h$ in (\ref{Hilbert1}) is regular too and its induced Finsler metric 
$F_h$  is determined  by
\be
F_h (y):= {1\over 2} \Big ( F_f(y)+F_f(-y) \Big ). \label{Hilbert2}
\ee 
T. Okada \cite{Ok} proved that the Funk metric $F_f$ satisfies the following equation
\be
{\pa F_f\over \pa x^i} = F_f {\pa F_f\over \pa y^i}.\label{Funkequation}
\ee
Okada uses  to prove the fact that  $F_f$ is of constant curvature $\kappa = -{1\over 4}$ and $F_h$ is of constant curvature $\kappa = -1$.
\section{Minkowski Spaces}

Minkowski spaces are finite dimensional vector spaces equipped with a Finsler metric 
invariant under translations. 
Thus Minkowski spaces are just vector spaces equipped with  Minkowski functionals.
 For a general  Finsler space $(M, F)$,
 each tangent space $T_xM$ with  $F_x:=F|_{T_xM}$ is a Minkowski space. Thus to study the geometric structure of a Finsler space, 
we need to study  Minkowski spaces first. 

Let $(V, F)$ be an $n$-dimensional  Minkowski space.
For each $y\in V \setminus \{0\}$, $F$ induces an inner product $g_y$ by 
(\ref{gy}).
$g_y$ satisfies the following homogeneity condition
\[
g_{\lambda y}(u, v) = g_y (u, v), \ \ \ \ \ \lambda >0.
\]
 Note that $g_y$ is independent of $y$ if and only if $F$ is Euclidean.
It  
 is natural to introduce 
the following quantity:
\be
{\bf C}_y (u, v, w):={1\over 2} {d \over dt} \Big [ g_{y+tw} (u, v) \Big ]\big  |_{t=0}.
\ee
The family  ${\bf C}:= \{ {\bf C}_y \}_{y\in V \setminus\{0\} }$ is called the
{\it Cartan torsion}. One can easily verify that ${\bf C}_y$ is a symmetric multi-linear form on $V$.
Moreover, ${\bf C}_y$ satisfies the following homogeneity
condition
\[
{\bf C}_{\lambda y} (u, v, w) = \lambda^{-1} {\bf C}_y (u, v, w), \ \ \ \ \ \lambda >0.\]
Note that  ${\bf C}=0$
if and only if $F$ is Euclidean.  Differentiating ${\bf C}_y$ with respect to $y$ yields 
a new quantity:
\be
\tilde{\bf C}_{y}(u, v, w, z):= 
{d \over dt} \Big [ {\bf C}_{y+tz}(u, v, w) \Big ]\big |_{t=0}.
\ee
Let $\tilde{\bf C}:=\{{\bf C}_y\}_{y\in V\setminus\{0\}}$. 
 $\tilde{\bf C}$ also
 gives us some geometric information on the Finsler metric \cite{Sh1}.

The mean of ${\bf C}_y$ is defined by
\be
{\bf I}_y (u):=\sum_{ij=1}^n g^{ij}(y) {\bf C}_y (e_i, e_j, u),
\ee
where $g_{ij}(y)=g_y(e_i, e_j)$. The family ${\bf I}=\{ {\bf I}_y \}_{y\in V\setminus\{0\}}$ is called
the {\it mean Cartan torsion}.
Deicke's Theorem \cite{De} says that 
${\bf C} =0$ if and only if ${\bf I}=0$.
Note that in dimension two, the family 
${\bf I}=\{{\bf I}_y \}_{y\in V\setminus\{0\}}$ completely determines the Cartan torsion.

\bigskip

There is another interesting quantity for Minkowski spaces associated with a Haar measure. 
Let $\mu$ be a Haar measure on $V$ which is invariant under translations.
Take an arbitrary  basis $\{e_i\}_{i=1}^n$ for $V$ and its dual basis $\{\omega^i\}_{i=1}^n $ for $V^*$,
$\mu$  can be expressed   by
$ d\mu = \sigma \; \omega^1 \wedge \cdots \wedge \omega^n. $
 We define
\be
\tau (y):= \ln {\sqrt{ \det\big ( g_{ij}(y) \big ) } \over \sigma },
\ee
where $g_{ij}(y):= g_y(e_i, e_j)$.
$\tau$  is a well-defined quantity
which is called the {\it distortion} of $(F, \mu)$ \cite{Sh2}\cite{Sh3}.
$\tau(y)$ satisfies the following homogeneity condition
\[ \tau(\lambda y) = \tau(y), \ \ \ \ \ \lambda >0.
\]
In general, $\tau(y)$ depends on the direction $y$.
Differentiating $\tau(y)$ with respect to $y$ yields the mean Cartan torsion. 
\be
{d \over dt} \Big [ \tau(y+tv) \Big ]\big |_{t=0} 
=  {\bf I}_y(v).
\ee
Therefore, the following  conditions are equivalent
(a) $\tau(y) = constant$;
(b) ${\bf I}=0$;
(c) ${\bf C} =0$;
(d) $F$ is Euclidean.

There are several special  Haar measures 
on a Minkowski space $(V, F)$. One of the natural Haar measures is the 
 {\it Busemann-Hausdorff}  measure $\mu_F$. $\mu_F$ can be expressed by
 $d\mu_F = \sigma_F \omega^1\wedge\cdots \wedge  \omega^n$, where
\be
\sigma_F := { {\rm Vol}({\rm B}^n) \over {\rm Vol} \{ (y^i) \in {\rm R}^n, \;
F(y^ie_i) < 1 \} } ,
\ee
where 
${\rm B}^n$ denote the unit ball in ${\rm R}^n$ and
${\rm Vol}$ denotes the Euclidean  measure on ${\rm R}^n$.
The Busemann-Hausdorff measure is the unique Haar measure
$\mu$ such that the unit ball ${\rm B}$ in $(V, F)$  has the same  volume as the standard unit ball
${\rm B}^n $ in ${\rm R}^n$.  It is somewhat surprising  that the Busemann-Hausdorff volume of
the Funk metric $F_f$ in (\ref{Funk2})  is finite. More precisely, for any metric $r$-ball $B(x, r)$ in the Funk space $(\Omega, F_f)$,
\[
\mu_F (B(x, r))= n \cdot 2^n\cdot {\rm Vol}({\rm B}^n) \int_0^{r/2} e^{-(n+1)t}\;  \sinh^{n-1}(t) dt \to {\rm Vol}({\rm B}^n). 
\]

\bigskip
Let $(V, F)$ be an $n$-dimensional Minkowski space and 
${\rm S}=F^{-1}(1)$ the indicatrix. 
There are  two induced metric structures  on the indicatrix ${\rm S}$. One is the Riemannian metric
$\dot{g}$ induced by $g_y$, and the other is the Finsler metric $\dot{F}$ induced by
$F$.

In 1949, L.A. Santal\'{o} \cite{Sa} proved that  if $F$ is reversible, then the Riemannian volume of the indicatrix ${\rm S}$ satisfies
\be
\mu_{\dot{g}} ({\rm S}) \leq {\rm Vol} ({\rm S}^{n-1} ),
\ee
equality holds if and only if $F$ is Euclidean.
However, there is no uniform lower bound on $ \mu_{\dot{g}}({\rm S})$.
For the further study on the Minkowski functional  $F$, one has to study
the geometry of $({\rm S}, \dot{g})$. It is surprising that the 
Riemannian curvature tensor  $\dot{\bf R}_y $ 
of $\dot{g}$ at $y\in {\rm S}$  takes a special form as follows:
\be
 \dot{\bf R}_y(u,v)w = {\bf C}_y( {\bf C}_y(u, w),v) -{\bf C}_y ({\bf C}_y(v, w), u)
+ \dot{g}_y (v, w) u - \dot{g}_y(u, w) v, 
\ee
where $u, v, w\in T_y{\rm S}\subset V$ and 
 ${\bf C}_y(u, v)=\xi$ is determined by $g_y(\xi, w): = {\bf C}_y(u, v, w)$.
The Brickell theorem says that in dimension  $n=\dim V > 2$,  
$\dot{g}$ has constant curvature $\kappa=1$ if and only if $F$ is Euclidean \cite{Bri}.

 For the Busemann-Hausdorff measure $\mu_{\dot{F}}$  on the indicatrix  ${\rm S}$, we have
 \be
 c_n \leq \mu_{\dot{F}} ({\rm S}) \leq c_n',   \label{vbh}
 \ee
 where $c_n$ and $c_n'$ are positive constants depending only on $n$. No sharp constants have 
 been determined in higher dimension. If $F$ is non-reversible, however, there is no
 uniform upper bound on $\mu_{\dot{F}}({\rm S})$.
For further investigation on the Minkowski functional  $F$, one has to study the
geometry of  $({\rm S}, \dot{F})$.  
Suppose that $\dot{F}$ is of constant curvature
$\kappa=1$. Is $F$ Euclidean ? 

\section{Connection and Geodesics}

Now we consider  general Finsler spaces. Geodesics are  the first objects coming to a geometer's sight
when he walks into an inner metric space.
 By definition, geodesics
are locally length-minimizing constant speed curves which are characterized locally by a system of second order ordinary differential equations.

Let $(M, F)$ be a Finsler space. For a $C^1$ curve $c: [a, b]\to M$,
the length of $c$ is given by
\[ \ell(c) =\int_a^b  F(\dot{c}(t)) dt.\]
A direct computation yields the Euler-Lagrange equations for 
a geodesic $c(t)$
\be
{d^2 x^i \over dt^2} + 2 G^i (\dot{c} ) =0,
\ee
where $(x^i(t))$ denote the coordinates of $c(t)$ and
$G^i$ in the standard local coordinate system $(x^i, y^i)$ in $TM$ are given by
\be
G^i(y) := {1\over 4} g^{il}(y) \Big \{ 
2 {\pa g_{jl}\over \pa x^k}(y) - {\pa g_{jk}\over \pa x^l} (y) \Big \} y^j y^k.
\label{Gi}
\ee 
where $g_{ij}(y)= g_y ({\pa \over \pa x^i}, {\pa \over \pa x^j})$.

A Finsler metric is said to be {\it positively complete} (resp. {\it complete}) if 
every geodesic on $(a, b)$  can be extended to a geodesic defined on
 $(a, \infty)$ (resp. $(-\infty, \infty)$). 
 The Funk metric in (\ref{Funk2}) is positively complete, but not complete,
 while the Hilbert metric in (\ref{Hilbert2}) is complete. Finsler metrics on 
 a compact manifold are always complete regardless the reversibility.
 
 \bigskip
With the  geodesic coefficients $G^i$ in (\ref{Gi}), we  define
 a map
 ${\rm D}_y : C^{\infty}(TM) \to T_xM $
for each $y\in T_xM$  by
 $$ {\rm D}_y U := \Big \{ dU^i (y) + U^j(x) {\pa G^i\over \pa y^j}(y) \Big \}
 {\pa \over \pa x^i}|_x,$$
where $U = U^i {\pa \over \pa x^i}\in C^{\infty}(TM)$.
${\rm D}_yU$ is called the {\it covariant derivative} of $U$ in the direction $y$.
We call the family ${\rm D}:= \{ {\rm D}_y \}_{y\in TM}$ the {\it canonical connection} of $F$.
W. Barthel first noticed this canonical connection.
 With this connection ${\rm D}$, we can define 
 the covariant derivative ${\rm D}_{\dot{c}} U(t)$  of a vector 
 field $U(t)$ along a curve $c(t)$, $ a \leq t \leq b$.
$U(t)$ is said to be {\it parallel} along $c$ if 
${\rm D}_{\dot{c}}U(t)=0$. Clearly, a curve $c$ is a geodesic if and only if 
the tangent vector field $\dot{c}(t)$ is parallel along $c$.
 The parallel translation 
 $P_c: T_{c(a)} M \to T_{c(b)}M$ is defined by
 \[ P( U(a)) = U(b)\]
 where $U(t)$ is  parallel along $c$. From the definition, we see that 
 $P_c$ is a linear transformation preserving  the inner products $g_{\dot{c}}$. In general, $P_c$ does not
 preserve the Minkowski functionals. We will discuss this issue in the next section.
 
It is natural to study the holonomy group defined by the above parallel
translations.
 A natural question is whether or not
 there are  more types of holonomy groups of Finsler spaces
  than the Riemannian 
 case. This problem remains open so far.

 \bigskip
 
 Let $(M, F)$ be a positively complete Finsler space. At each point
 $x\in M$, we define
 a map $\exp_x: T_xM \to M$ by
 \[ \exp_x (y) := c(1),\]
 where $c(t)$ is the geodesic with $\dot{c}(0)=y$. The Hopf-Rinow theorem says that
 $\exp_x$ is onto for all $x\in M$. $\exp_x$ is called the {\it exponential map}
 at $x$. From the O.D.E. theory, J.H.C. Whitehead \cite{Wh}
proved  that $\exp_x$ is $C^{\infty}$ on $T_xM\setminus\{0\}$ and only $C^1$ at 
 the origin. Akbar-Zadeh \cite{AZ}  proved that $\exp_x$ is $C^2$ at the origin for all $x$ if and only if
 ${\rm D}$ is an  affine connection.

 \bigskip\section{Non-Riemannian Curvatures}

 The canonical connection ${\rm D}$  has all the properties of an affine connection
 except for the linearity in $y$. Namely,
 ${\rm D}_{y_1+y_2} \not={\rm D}_{y_1}+{\rm D}_{y_2}$ in general. 
 To measure the non-linearity, it is natural to introduce
 the following quantity \cite{Sh1}
 \be
 {\bf B}_y (u, v, w) := {\pa^2 \over \pa s \pa t } \Big [
 {\rm  D}_{y+sv+tw}U \Big ]|_{s=t=0},
 \ee
 where $U\in C^{\infty}(TM)$ with $U(x)=u$.
 One can easily  verify that ${\bf B}_y$ is a symmetric multi-linear form
 on $T_xM$. We call the family
 ${\bf B} := \{ {\bf B}_y \}_{y\in TM \setminus\{0\}}$ the {\it Berwald curvature}.
 A Finsler metric is called a {\it Berwald metric} if ${\bf B}=0$.
 L. Berwald proved a simple fact that  ${\bf B}=0$ if and only if 
  ${\rm D}$ is an affine connection.

 For Riemannian metrics, ${\bf B} =0$ and ${\rm D}$ is just the Levi-Civita 
 connection.
 There are non-Riemannian Berwald metrics with ${\bf B}=0$.
 Consider the following type of Finsler metric:
 \be
 F(y):= \alpha(y) + \beta(y), \label{Randers}
 \ee
where $\alpha(y):=\sqrt{g(y, y)}$ is a Riemannian metric and 
$\beta(y)$ is a $1$-form with 
$\alpha$-length $\|\beta\|<1$. $F$ is called a {\it Randers metric}.
M Hashiguchi and Y. Ichijy\={o} \cite{HaIc1} first noticed that 
if $\beta$ is parallel with respect to $\alpha$, then $F=\alpha +\beta$
is a Berwald metric.
Later,  they   proved that if $d\beta =0$, then
 $F=\alpha+\beta$ has the same geodesics as  $\alpha$
and vice versa \cite{HaIc2}. 
 
 Y. Ichij\={o} \cite{Ic}  proved that on  a Berwald space, 
 the parallel translation along any geodesic preserves the Minkowski
 functionals. Thus Berwald spaces can be viewed as 
 Finsler spaces modeled on a single Minkowski space. According to Szab\'{o} \cite{Sz}, if a Finsler metric 
 $F$ is Berwaldian, then there is a Riemannian metric $g$ 
 whose Levi-Civita connection coincides 
 with the canonical connection  of $F$.

 Define the mean of ${\bf B}_y$  by
 \be
 {\bf E}_y (u, v) := {1\over 2}
 \sum_{i=1}^n g^{ij}(y) g_y \Big ( {\bf B}_y (u, v, e_i ) , e_j \Big ),
 \ee
 where $g_{ij}(y)=g_y(e_i, e_j)$. 
The family ${\bf E}=\{ {\bf E}_y \}_{y\in TM\setminus\{0\}}$ is called the {\it mean Berwald curvature}. 
${\bf E}$ is also related to the S-curvature ${\bf S}$. See \cite{Sh1}  and 
(\ref{ES}) below.

 \bigskip
 As we have mentioned above, the parallel translation along curve in a Berwald space 
 preserves the Minkowski functionals. Thus the Cartan torsion in a Berwald space
 does not change along geodesics. 
 To measure the rate of changes of the Cartan torsion
along geodesics in a general Finsler space, we will introduce a weaker quantity
than the Berwald curvature. 
 For a vector $y\in T_xM\setminus\{0\}$, let
 $c(t)$ denote the geodesic with $\dot{c}(0)=y$.
 Take arbitrary vectors $u, v, w\in T_xM$ and extend them to
 parallel vector fields $U(t), V(t), W(t)$ along $c$.
 Define
 \be
 {\bf L}_y (u, v, w):=  {d\over dt}
 \Big [ {\bf C}_{\dot{c}(t)}\Big (U(t), V(t), W(t)\Big )\Big ]|_{t=0}.\label{LC}
 \ee
 The family ${\bf L}:= \{ {\bf L}_y \}_{y\in TM\setminus\{0\}}$
 is called the {\it Landsberg curvature}. A Finsler metric 
is called a {\it Landsberg metric} if ${\bf L}=0$ \cite{Sh1}. Landsberg metrics form 
an important class of Finsler spaces.
 We have the following equation \cite{Sh1}
 \be
 {\bf L}_y (u, v, w) = - {1\over 2} g_y \Big ( {\bf B}_y (u, v, w), y \Big ).\label{JB}
 \ee
 From (\ref{JB}), we immediately conclude that
 every Berwald space is a Landsberg space.
 It is an open problem in Finsler geometry 
 whether or not there is a Landsberg metric which is not a Berwald metric. 
 So far no example has been found.
Differentiating  ${\bf L}$ along geodesics yields a new quantity:
\be
\dot{\bf L}_y (u, v, w):= {d\over dt} \Big [
{\bf L}_{\dot{c}(t)} \Big (U(t), V(t), W(t)\Big )\Big ]_{t=0}.\label{dotJ}
\ee

 \bigskip
Using (\ref{Funkequation}), we can show that
the Funk metric  $F=F_f$  in (\ref{Funk2}) satisfies
\be
{\bf L}_y(u,v,w) + {1\over 2} F(y) {\bf C}_y(u,v,w) =0,\label{ll}
\ee
and  the Hilbert metric in (\ref{Hilbert2}) satisfies  
\be
\dot{\bf L}_y(u, v, w) - F^2 (y) {\bf C}_y(u,v,w) =0. \label{kk}
\ee

 \bigskip
 The Landsberg curvature
 ${\bf L}_y$ satisfies the following homogeneity condition
 \be
 {\bf L}_{\lambda y}(u,v,w) = {\bf L}_y(u,v,w), \ \ \ \ \ \lambda >0.
 \ee
 In general,  ${\bf L}_y$ depends on the direction $y$. 
Differentiating ${\bf L}_y$ with respect to $y$ yields another quantity \cite{Sh1}
 \be
 \tilde{\bf L}_y(u, v, w, z):=
 {d \over dt} \Big [ {\bf L}_{y+tz} (u, v, w) \Big ]|_{t=0}.
\ee
 One can easily verify that $\tilde{\bf L}=0$ if and only if  
 ${\bf L}=0$. When ${\bf L}\not=0$, 
$\tilde{\bf L}$ gives us some other geometric information on the Finsler metric.

 \bigskip

 Define the mean of ${\bf L}_y$  by
 \be
 {\bf J}_y (u):= \sum_{ij=1}^n g^{ij}(y) {\bf L}_y (u, e_i, e_j ).
 \ee
 The family ${\bf J}=\{{\bf J}_y\}_{y\in TM\setminus\{0\}}$ 
 is called the {\it mean Landsberg curvature} \cite{Sh1}. 
 From the definitions of ${\bf I}$
and ${\bf J}$, we have
\be
{\bf J}_y (u) =  {d\over dt} \Big [ {\bf I}_{\dot{c}(t)}\Big (U(t)\Big ) \Big ]\big|_{t=0}
,
\ee
where
 $c(t)$ is the geodesic 
with $\dot{c}(0)=y$ and $U(t)$ is a parallel vector field along $c$ with $U(0)=u$.
In dimension two, ${\bf J}$ completely determines ${\bf L}$.  It is an interesting 
problem to study the difference between Finsler metrics with
${\bf J}=0$ and those with ${\bf L}=0$.
\bigskip

 There is an induced Riemannian metric  of Sasaki type on $TM\setminus\{0\}$.
T. Aikou proved that if  ${\bf L}=0$, then 
  all the slit tangent spaces $T_xM\setminus\{0\}$ are totally geodesic
 in $TM\setminus\{0\}$ \cite{Ai}. Along the same line, one can show that if 
 ${\bf J}=0$, then  all the slit tangent spaces $T_xM\setminus\{0\}$ are minimal
 in $TM\setminus\{0\}$.

\bigskip
Consider an arbitrary regular measure $\mu$ on a Finsler 
space $(M, F)$. $\mu$ induces a Haar measure $\mu_x$ in each tangent space $T_xM$. 
Hence the distortion $\tau$ is defined for $(T_xM, F_x, \mu_x)$. 
To measure the rate of changes of the distortion along geodesics, we define
 \be
 {\bf S}(y):= {d \over dt} \Big [ \tau \Big (\dot{c}(t)\Big ) \Big ]_{t=0},
 \ee
where $c(t)$ is the geodesic with $\dot{c}(0)=y$.
We call the scalar function ${\bf S}$ the {\it S-curvature} \cite{Sh1}\cite{Sh3}. 
Differentiating the S-curvature along geodesics yields a new quantity:
\be
\dot{\bf S}(y):= 
{d\over dt} \Big [ {\bf S}(\dot{c}(t)) \Big ]\big |_{t=0}.
\ee
See \cite{Sh2} for further discussions.
When ${\bf S}\not=0$,  $\dot{\bf S}$ gives us some other  geometric information on the Finsler metric $F$ and the regular measure $\mu$.
See (\ref{volumeTaylor}) below.

\bigskip

The S-curvature ${\bf S}(y)$ satisfies the following homogeneity condition
 \be
 {\bf S}(\lambda y ) = \lambda {\bf S}(y), \ \ \ \ \lambda >0.
 \ee
 In general, ${\bf S}(y)$ is not linear in $y$.
 Differentiating it twice with respect to $y$ gives  no new quantity. Namely, we have 
\be
{\bf E}_y (u, v) = {1\over 2} {\pa^2 \over \pa s \pa t} \Big [ {\bf S}(y+ su+tv) \Big ]_{s=t=0} . \label{ES}
 \ee
 Thus ${\bf S}(y)$ is linear in $y\in T_xM$ if and only if ${\bf E}=0$ on $ T_xM\setminus\{0\}$. In particular, if $F$ is a
Berwald metric, then ${\bf S}(y)$ is linear 
in $y\in T_xM$ for all $x$ \cite{Sh1}. 
In fact,  ${\bf S}=0$ for Berwald metrics
if we consider the S-curvature of  the Busemann-Hausdorff
measure $\mu_F$. This fact is proved by the author \cite{Sh3}.   
Finsler spaces with ${\bf E}=0$ deserve further 
investigation. There are some non-Berwaldian Randers metrics
 with ${\bf E}=0$ and ${\bf S}=0$. 
 For the Funk metric $F=F_f$ in (\ref{Funk2}), the S-curvature and the mean Berwald curvature are constant in the following sense.
\begin{eqnarray*}
&&{\bf S}(y)= {n+1\over 2} F(y),\\
&& {\bf E}_y(u,v)  = {n+1\over 4 F^3(y)} \Big \{F^2(y) g_y (u, v) - g_y(y, u)g_y(y, v)\Big \} .
\end{eqnarray*}
This is proved in \cite{Sh1}.

\section{Riemann Curvature}

As matter of fact, all the quantities defined in the 
previous sections vanish on a Riemannian space. Thus we do not see
these non-Riemannian quantities  at all  in Riemannian geometry. 
 A. Einstein used Riemannian geometry to describe his
general relativity theory, assuming that a spacetime is {\it always} Riemannian.

 For Riemannian spaces,
 there is only one notion of curvature---Riemann curvature, that was introduced by 
B. Riemann in 1854 as a generalization of the Gauss curvature for surfaces. Since then, the Riemann curvature
 became the central concept in Riemannian geometry.
 Due to the efforts by L. Berwald in 1920's,
 the Riemann curvature can be extended
to the Finslerian case \cite{Ber}.

Let 
$(M, g)$ be a Riemannian space and ${\rm D}$ denote the Levi-Civita connection of $g$. 
The {\it Riemann curvature tensor }
is defined by
\[ {\rm R}(u, v)w :=\Big \{ {\rm D}_U {\rm D}_V W 
- {\rm D}_V {\rm D}_U W - {\rm D}_{[U, V]} W \Big \}|_x,
\ \ \ \ u,v,w\in T_xM,
\]
where $U,V, W$ are  local vector fields with $U(x)=u, V(x)=v, W(x)=w$. 
The core part of the Riemann curvature tensor is the following quantity:
\be
{\bf R}_y (u):= {\rm R}(u, y)y.
\ee
The Riemann curvature ${\bf R}_y: T_xM \to T_xM$ is a self-adjoint linear
transformation with respect to $g$ and it satisfies  ${\bf R}_y(y)=0$.
The family ${\bf R}=\{ {\bf R}_y \}_{y\in TM\setminus\{0\}}$ 
is called the {\it Riemann curvature}. 
With a little trick by the author, one  can extend the notion of 
Riemann curvature to  Finsler metrics without employing connections on the slit
tangent bundle $TM\setminus\{0\}$. 

Let $(M, F)$ be a Finsler space. Given a vector $y\in T_xM \setminus\{0\}$, extend it to a local
nowhere zero {\it geodesic field} $Y$ (i.e., all integral curves of $Y$ are geodesics). 
$Y$ induces a Riemannian metric 
\[ \hat{g}:= g_Y.\]
Let $\hat{\bf R}$ denote the Riemann curvature of $\hat{g}$ as defined above. Define
\be
{\bf R}_y := \hat{\bf R}_y .
\ee
One can verify that ${\bf R}_y$ is independent of the geodesic
extension $Y$ of $y$. Moreover,
${\bf R}_y$ is self-adjoint with respect to
$g_y$, i.e., 
\[ g_y \Big ({\bf R}_y(u), v\Big ) = g_y \Big (u, {\bf R}_y (v)\Big ),\]
and it satisfies  ${\bf R}_y(y)=0$ \cite{Sh1}\cite{Sh2}. Let
$W_y:= \{ u\in T_xM, \; g_y(y, u)=0 \}$.
Then
${\bf R}_y|_{W_y}: W_y \to W_y$
is again a self-adjoint linear transformation 
with respect to $g_y$. Denote the eigenvalues  
 of ${\bf R}_y|_{W_y}$ by
 \[\kappa_1 (y)
\leq \cdots \leq \kappa_{n-1}(y).\]
 They   are the most important intrinsic invariants of 
the Finsler metric. We call $\kappa_i(y)$  the $i$-th
principal curvature in the direction $y$. The trace of 
${\bf R}_y$ is denoted by ${\bf Ric}(y)$  which is called the {\it Ricci curvature}.
${\bf Ric}(y)$ is given by
\be
{\bf Ric}(y) :=\sum_{ij=1}^n  g^{ij}(y) 
g_y\Big  ( {\bf R}_y (e_i), e_j\Big ) = \sum_{i=1}^{n-1} \kappa_i (y).
\ee

\bigskip
The Ricci curvature and the S-curvature determine the local behavior
of the Busemann-Hausdorff measure of small  metric balls around a point.
Let ${\rm B}_x$ denote the unit ball in $(T_xM, F_x)$ and $\mu_x$ the 
induced Busemann-Hausdorff measure of $F_x$ on $T_xM$.
 Assume that $F$ is reversible. Then the
Taylor expansion of $\mu_F(B(x, \e))$ of a small
metric ball $B(x, \e)$ is given by
\be
\mu_F (B(x, \e))
= {\rm Vol}({\rm B}^n)
\Big \{ 1 - {1\over 6(n+2)} r(x) \; \e^2 + O(\e^3) \Big \}, \label{volumeTaylor}
\ee
where
\be
r(x) := {n+2 \over n\cdot {\rm Vol}({\rm B}^n)}\int_{{\rm B}_x} \Big \{
{\bf Ric}(y) d\mu_x    + 3n
 \Big [ \dot{\bf S}(y) - {\bf S}^2 (y) \Big ]\Big \}
d\mu_x .
\ee
See \cite{Sh2} for details.

\section{Constant Curvature}

Now let  us take a close look at   Finsler spaces   of constant curvature $\kappa$.
A Finsler metric is said to be {\it of scalar curvature} 
if there is a scalar function $\kappa(y)$ on
$TM \setminus \{0\}$ such that for any $y\in T_xM \setminus\{0\}$, the principal 
curvatures $\kappa_i(y)=\kappa(y)$, $i=1, \cdots, n-1$. 
By definition, all two dimensional Finsler metrics
are of scalar curvature $\kappa(y)$. $F$ is said to be of {\it constant
curvature} $\kappa$ (resp. {\it  constant Ricci curvature}) if  $\kappa_i(y) =\kappa, \; i =1, \cdots, n-1$ (resp. $\sum_{i=1}^{n-1} \kappa_i(y) = (n-1)\kappa$).

We have the following important 
equation \cite{AZ}
\be
\dot{\bf L}_y (u,v,w)+ \kappa\;  F^2(y) {\bf C}_y (u,v,w)=0. \label{dotLC}
\label{Numata}
\ee
The Cartan torsion and Landsberg curvature take special values along 
geodesics.
Let $c(t)$ be 
an arbitrary {\it unit speed}  curve. Take a  parallel vector field $V(t)$  along 
$c(t)$.
Let
\be
{\bf C}(t) := {\bf C}_{\dot{c}(t)} (V(t), V(t), V(t)).
\ee
From (\ref{dotJ}) and (\ref{dotLC}), we obtain
the following 
important equation \cite{Nu}\cite{AZ}
\be
{\bf C}''(t) + \kappa {\bf C}(t)=0. \label{CCC}
\ee
This immediately implies that Landsberg space of 
constant curvature $\kappa \not=0$ must be Riemannian. This is observed by S. Numata \cite{Nu}.
Solving (\ref{CCC}), we obtain
\be
{\bf C}(t) =\cases{a \sinh(t) + b \cosh(t), & if $\kappa =-1$, \cr\\
at + b, & if $\kappa =0$,\cr\\
a \sin(t) + b \cos (t), & if $\kappa = 1$.}\label{C(t)}
\ee
Define ${\bf L}(t)$ in the same way as above for the Landsberg
curvature. From the definition of ${\bf L}$,
we have
${\bf L}(t) = {\bf C}'(t)$. Then we obtain a
formula for ${\bf L}(t)$ \cite{AZ}.

\bigskip
Take two parallel vector fields $V(t)$ and $W(t)$ along $c$. Assume that
both $V(t)$ and $W(t)$ 
are $g_{\dot{c}(t)}$-orthogonal to $\dot{c}(t)$
for some $t=t_o$ (hence for all $t$). Set
\[
\tilde{\bf C}(t) : =\tilde{\bf C}_{\dot{c}(t)} (V(t), V(t), V(t), W(t)).
\]
By studying the Ricci identities
and the Bianchi identities, we
obtain
\be
\tilde{\bf C}(t) =
\cases{ a \sinh(2t) + b\cosh(2t) + c, & if $\kappa =-1$, \cr\\
at^2 + bt +c , & if $\kappa=0$, \cr\\
a\sin (2t) + b \cos (2t) +c,& if $\kappa =1$.}\label{tildeC(t)}
\ee
Define $\tilde{\bf L}(t)$ in the same way as above  for $\tilde{\bf L}$ . We can show that
$ \tilde{\bf L}(t) = \tilde{\bf C}'(t) + c'$. Then we 
obtain a formula for $\tilde{\bf L}(t)$ \cite{Sh1}.

\bigskip

Complete Finsler metrics of constant curvature $\kappa <0$ must be Riemannian if the Cartan torsion 
does not grow exponentially. This fact is due to
Akbar-Zadeh \cite{AZ}. Using (\ref{Funkequation}), T. Okada \cite{Ok}
verified that 
the Funk metric $F_f$ in (\ref{Funk2})  is of constant curvature
$\kappa =  - {1\over 4}$ and
the Hilbert metric $F_h$ in (\ref{Hilbert2}) is of constant curvature
$\kappa = -1$. 
By (\ref{ll}), we can show that 
 the  Cartan torsion of $F_f$ is bounded  along any geodesic. Note that $F_f$ is not Riemannian !
because it is only positively complete. 
 Since $F_h$ is non-Riemannian,  
the Cartan torsion of $F_h$ must grow exponentially along geodesics in one direction.

\bigskip
Positively complete Finsler spaces of constant curvature $\kappa=0$ must be locally Minkowski
if ${\bf C}$ and $\tilde{\bf C}$ are bounded along geodesics. 
This fact is also due to Akbar-Zadeh \cite{AZ}. So far, we do not know if there are any 
positively complete
Finsler spaces of constant curvature $\kappa=0$, except for locally Minkowski spaces.

\bigskip

There are infinitely many projectively flat 
Finsler metrics of constant curvature $\kappa = 1$ on ${\rm S}^n$
constructed by  R. Bryant \cite{Br1}\cite{Br2} recently. Bryant metrics are non-reversible. 
So far, no reversible Finsler metric of constant curvature 
$\kappa =1$ has been found on ${\rm S}^n$, except for the standard Riemannian metric. 
The author can prove that for  any Finsler metric  on a simply connected compact manifold $M$,
if it has constant curvature $\kappa = 1$, then
$M$ must be diffeomorphic to ${\rm S}^n$ and geodesics
are all closed with length of $2\pi$. From (\ref{C(t)})
and (\ref{tildeC(t)}), we see that 
${\bf C}$ has period of $2\pi$ on parallel vector fields
along any 
geodesic, while $\tilde{\bf C}$
has period of $\pi$ on parallel vector fields
orthogonal to the geodesic \cite{Sh1}.

\bigskip 

All known Finsler metrics of constant curvature are locally  projectively flat, i.e.,
at every point, there is a local  coordinate system in which  the geodesics are 
straight lines. It is an interesting problem to find
Finsler metrics of constant curvature without this property.

Consider two  pointwise projectively related Finsler metrics $F$ and $\tilde{F}$ on a manifold.
Suppose that $F$ and $\tilde{F}$ has constant Ricci
curvature $\kappa$ and $\tilde{\kappa}$, respectively.
Then using A. Rapcs\'{a}k's  equation, we can show that  for any unit speed geodesic $c(t)$ of $F$,
the function $\varphi(t):=1/\sqrt{\tilde{F}(\dot{c}(t))}$ satisfies
\be
\varphi''(t) + \kappa \varphi (t) = {\tilde{\kappa}\over \varphi^3(t)}. \label{Einstein}
\ee
See \cite{Sh1}\cite{Sh4}. 
By (\ref{Einstein}), we can show that
the Hilbert metric is the only
complete, reversible, projectively flat 
Finsler  metric of constant curvature $\kappa =-1$
on a strongly convex domain in ${\rm R}^n$. There are might be many
positively complete non-reversible  projectively flat Finsler metrics of constant curvature
$\kappa = -{1\over 4}$ on a strongly convex domain in ${\rm R}^n$.
So far we only have the Funk metric with this property.

It is an open problem whether or not there is a (positively) 
complete Finsler space which does not admit 
any (positively) complete Finsler metrics of scalar curvature.  This leads to the study on the topology
of (positively) complete Finsler spaces of scalar curvature.

\section{Comparison Geometry}
In this section, we will discuss several global results using 
comparison techniques. 

Let $(M, F)$ be a positively complete Finsler space.
Take a geodesic variation $c_s (t)$ of a geodesic $c(t)$, i.e.,
$c_0(t)=c(t)$ and each $c_s(t)$ is a geodesic. Let
$J(t):= {\pa c_s \over \pa s}|_{s=0}(t)$. $J(t)$ is a vector field 
  along $c$ which is called a {\it Jacobi field}. The behavior of $J(t)$
along $c$ is controlled by the following Jacobi  equation
\be
{\rm D}_{\dot{c}}{\rm D}_{\dot{c}} J(t) + {\bf R}_{\dot{c}(t)} (J(t)) =0.
\ee

Take a geodesic $c(t)=\exp_x(ty), 0 \leq t < \infty$ and a special
geodesic variation $c_s (t):= \exp_x (t (y+sv))$.  
The standard comparison argument by Cartan-Hadamard and Bonnet-Meyers
gives the following
important global results
in comparison Finsler geometry.

\begin{thm} {\rm (\cite{Aus})}  Let $(M, F)$ be a positively complete
Finsler space. Suppose that the Riemann curvature is nonpositive, i.e,
the principal curvatures 
\[
\kappa_i (y) \leq 0, \ \ \ \ \ \  i= 1, \cdots, n-1.\] 
Then the exponential map  $\exp_x: T_xM \to M$ is
an onto covering map.  Thus $M$ is a $K(\pi, 1)$ space.
\end{thm}

\begin{thm} {\rm (\cite{Aus})} Let $(M, F)$ be a positively
complete Finsler space. Suppose that 
the Ricci curvature is strictly positive., i.e.,
there is a positive constant $\lambda$ such that
\[ \sum_{i=1}^{n-1} \kappa_i(y)  \geq (n-1) \lambda >0.\]
Then the exponential map
$\exp_x: T_xM \to M$ is singular at $r y$
for any unit vector $y\in T_xM$ at $r \leq \pi/\sqrt{\lambda}$. Thus  the 
diameter of $M$ and its universal cover $\tilde{M}$ is bounded by
${\rm Diam}(M)\leq \pi/\sqrt{\lambda}$, and   the fundamental group $\pi_1(M)$ must be finite.
\end{thm}

Applying the Morse theory to  the loop space, one can  prove the following theorem for homotopy groups.
\begin{thm}   
Let $(M, F)$ be a compact simply connected Finsler space. Suppose that the principal 
curvature $ \kappa_1 (y) \leq \cdots\leq  \kappa_{n-1}(y) $ satisfy the following pinching condition
for some  $ 2 \leq k \leq n-2$, 
\be
 {1\over 4}  < {1\over k} \sum_{i=1}^k \kappa_i(y), \ \ \ \ \ \kappa_{n-1}(y) \leq 1.
 \ee
 Then $\pi_i(M) =0$ for $i =1, \cdots, n-k$.
 \end{thm} 

\bigskip

Let $(M, F)$ be a positively complete space. 
Define by $B(x, r)$ and $S(x, r)$ the metric ball and sphere around $x$ with radius $r$, respectively.
There is a naturally induced measure $\nu_F$ on the regular part of $S(x, r)$ such that the coarea formula holds
\be
 \mu_F (B(x, r)) =\int_0^r \nu_F (S(x, t)) dt.\label{coarea}
\ee
Let $\mu_{\dot{F}}$ denote  the Busemann-Hausdorff measure of the induced 
Finsler metric $\dot{F}$ on $S(x, r)$. In general,
$\nu_F \not= \mu_{\dot{F}}$. If $F$ is reversible, then 
\be
 c_n\;  \mu_{\dot{F}}\leq \nu_F  \leq  c_n'\;  \mu_{\dot{F}},\label{vhv}
\ee
 where $c_n, c_n'$ are positive constants. 
If $F$ is non-reversible, the  inequality on the left side of (\ref{vhv}) does  not hold.
The coarea formula (\ref{coarea}) together with (\ref{vhv}) implies (\ref{vbh}). 
See \cite{Sh2} for more details.
Further estimates on the geometry of $S(x, r)$ give the following comparison
result on the Busemann-Hausdorff measure $\mu_F$ under certain curvature bounds.

\begin{thm}\label{thm7.4} {\rm (\cite{Sh2}\cite{Sh3})} Let $(M, F)$ be an $n$-dimensional  positively complete
Finsler space. Suppose that the Ricci curvature and
the S-curvature satisfy
\be
{\bf Ric}/F^2 \geq (n-1)\lambda,
\ \ \  \ {\bf S}/F \geq (n-1)\delta,\label{RicS}
\ee
where $\lambda, \delta$ are positive  constants. Then 
the ratios
$\mu_F (B(x, r))/V_{\lambda,\delta}(r)$ 
and $\nu_F (S(x, r))/V'_{\lambda,\delta}(r)$ are non-increasing, where
\[ V_{\lambda,\delta}(r):={\rm Vol}({\rm S}^{n-1}) 
\int_0^r \Big [ e^{-\delta t}\; {\bf s}_{\lambda}(t) \Big ] ^{n-1}  dt,\]
and ${\bf s}_{\lambda}(t)$ satisfies
\[ {\bf s}_{\lambda}''(t) +\lambda \; {\bf s}_{\lambda}(t)=0,
\ \ \ \ \ {\bf s}_{\lambda}(0)=0, \ {\bf s}'_{\lambda}(0)=1.\]
\end{thm}

Theorem \ref{thm7.4} has a number of applications.
Let $M$ be a compact oriented manifold. The 
canonical $L^1$-norm $\|\cdot \|_1$ on the complex 
$C_k(M)$ of singular real  chains is defined by
\[ \|c \|_1 :=\sum_i |r_i|, \ \ \ \ \ c=\sum_i r_i \sigma_i.\]
For a real homology class $z\in H_k(M)$, define
\[ \|z\|_1 = \inf_{z=[c]} \|c\|_1 .\]
For the fundamental class $[M]\in H_n(M)$, let
\[ \|M\|:= \|[M]\|_1.\]
 $\|M\|$ is called the {\it Gromov invariant} of $M$.
$\|M\|$ is not necessarily an integer. Gromov proved that if 
 $\pi_1(M)$ is amenable, then $\|M\| =0$. 

\begin{thm} Let $(M, F)$ be an $n$-dimensional  reversible
compact Finsler space. Suppose that the Ricci curvature
and the S-curvature satisfy
 the bounds (\ref{RicS}) with
$\lambda, \delta \leq 0$. Then 
\be
\| M\| \leq n! (n-1)^n (\sqrt{|\lambda|}+|\delta| )^n \mu_F(M).\label{min1}
\ee
Further, there is a constant $\e(n) >0$ if
\be
 (\sqrt{|\lambda|}+|\delta| )^n \mu_F(M) \leq \e(n),\label{min2}
 \ee
then $\|M\|=0$.
\end{thm}
The theorem for Riemannian spaces was proved by M. Gromov \cite{Gr}. The proof for the general case follows from
Gromov's argument by using Theorem \ref{thm7.4}.

\bigskip

\noindent
zshen@math.iupui.edu

\begin{thebibliography}{AnBnCng3}
\bibitem[Ai]{Ai} T. Aikou,
{\it Some remarks on the geometry of tangent bundles of Finsler spaces},
Tensor, N. S. {\bf 52}(1993), 234-242.

\bibitem[Aus]{Aus}
L. Auslander, {\it On curvature in
Finsler geometry}, Trans. Amer. Math. Soc. 
{\bf 79}(1955), 378-388.


\bibitem[AZ]{AZ}H. Akbar-Zadeh, 
{\it  Sur les espaces de Finsler \'{a} courbures sectionnelles constantes}, 
Bull. Acad. Roy. Bel. Cl, Sci, 5e S\'{e}rie - Tome LXXXIV 
(1988) 281-322.

\bibitem[Ber]{Ber}
L. Berwald, {\it Untersuchung der Kr\"{u}mmung allgemeiner metrischer R\"{a}ume auf Grund des in ihnen herrschenden Parallelismus}, Math. Z.
{\bf 25}(1926), 40-73.


\bibitem[Bri]{Bri} F. Brickell,
{\it A theorem on homogeneous functions}, J. London Math. Soc.
{\bf 42}(1967), 325-329.


\bibitem[Br1]{Br1}  {R. Bryant, \textit{Finsler structures on the
2-sphere satisfying $K=1$}, Finsler Geometry, Contemporary Mathematics {\bf 196},
Amer. Math. Soc., Providence, RI, 1996, 27-42. }

\bibitem[Br2]{Br2}  {R. Bryant, \textit{Projectively flat Finsler $2$%
-spheres of constant curvature}, Selecta Math., New Series, 
{\bf 3}(1997), 161-204. }


\bibitem[De]{De} A. Deicke, {\it \"{U}ber die Finsler-R\"{a}ume mit $A_i=0$}, Arch. Math. {\bf 4}(1953), 45-51.

\bibitem[Funk]{Funk} P. Funk, {\"{U}ber Geometrien, bei denen
die Geraden die K\"{u}rzesten sind}, Math. Ann. {\bf 101}(1929), 226-237.

\bibitem[Gr]{Gr} M. Gromov, {\it Volume and bounded cohomology}, I. H. E. S. Publ. Math. {\bf 56}(1983), 213-307.

\bibitem[HaIc1]{HaIc1} 
M. Hashiguchi and Y. Ichijy$\bar{o}$,
{\it On some special $(\alpha, \beta)$ metrics},
Rep. Fac. Sci. Kagoshima Univ. {\bf 8}(1975), 39-46.

\bibitem[HaIc2]{HaIc2} 
M. Hashiguchi and Y. Ichijy\={o},
{\it Randers spaces with rectilinear geodesics}, 
Rep. Fac. Sci. Kagoshima Univ. (Math. Phys. \& Chen.),
{\bf 13}(1980), 33-40. 

\bibitem[Ic]{Ic} 
 Y. Ichijy\={o},
{\it  Finsler spaces modeled on a Minkowski space}, 
J. Math. Kyoto Univ. {\bf 16}(1976), 
 639--652.

\bibitem[Nu]{Nu} S. Numata, {\it 
On Landsberg spaces of scalar curvature}, J. Korea Math. Soc. {\bf 12}(1975), 97-100. 

\bibitem[Ok]{Ok} T. Okada, {\it On models of projectively flat Finsler spaces of constant negative curvature}, Tensor, N. S. {\bf 40}(1983), 117-123.

\bibitem[Sa]{Sa} L.A. Santal\`{o}, {\it Un invariante afin para los cuerpos convexos del espacio de n dimensiones}, Portugaliae Mathematica,
{\bf 8} (1949), 154-161.

\bibitem[Sh1]{Sh1} Z. Shen, 
{\it Differenial Geometry of Spray and Finsler Spaces},
Kluwer Academic Publishers, 2001.

\bibitem[Sh2]{Sh2} Z. Shen,
{\it Lectures on Finsler Geometry},
in preparation.

\bibitem[Sh3]{Sh3} Z. Shen,
{\it Volume comparison  and its applications in Riemann-Finsler
geometry}, Advances in Math.
{\bf 128}(1997), 306-328.

\bibitem[Sh4]{Sh4} Z. Shen, {\it On projectively related Einstein metrics in Riemann-Finsler geometry},
to appear in Math. Ann.

\bibitem[Sz]{Sz} Z. Szab\'{o}, 
{\it Positive definite Berwald spaces (Structure theorems on Berwald spaces)}, 
 Tensor,  N. S.
{\bf  35}(1981),  25-39. 

\bibitem[Wh]{Wh}
J.H.C. Whitehead, {\it Convex regions in the geometry of paths}, Quart. J. Math. Oxford Ser. {\bf 3}(1932), 33-42.


\end{thebibliography}
\end{document}